\documentclass[11pt]{article}
\usepackage[a4paper]{geometry}

\usepackage{amsmath}
\usepackage{amssymb}
\usepackage{amsthm}
\usepackage{mathrsfs}
\usepackage{bbm}
\usepackage{empheq}
\usepackage[loose]{subfigure}
\usepackage{epsfig}
\usepackage{graphicx}
\usepackage{psfrag}
\usepackage[usenames,dvipsnames]{pstricks}
\usepackage{pst-plot}
\usepackage[colorlinks]{hyperref}
\usepackage{tabls}
\usepackage{paralist}
\usepackage{hyperref}

\DeclareSymbolFontAlphabet{\Bbb}{AMSb}

\newlength{\fixboxwidth}
\newcommand{\argmin}{\mathop{\mathrm{arg\,min}}}
\newcommand{\COMMENT}[1]{}
\newcommand{\E}{\mathbb{E}}
\newcommand{\Er}{\mathcal{E}}

\newcommand{\Dmaps}{\mathfrak{D}}
\newcommand{\Dmap}{\mathbb{D}}
\newcommand{\Drv}{D}
\newcommand{\Dspace}{\mathcal{D}}
\newcommand{\Ddata}{d}

\newcommand{\one}{\mathbbm{1}}
\renewcommand{\P}{\mathbb{P}}

\newcommand{\R}{\mathbb{R}}

\newcommand{\X}{\mathcal{X}}

\DeclareMathOperator{\Image}{Image}

\newcommand{\quark}{\setbox0\hbox{$x$}\hbox to\wd0{\hss$\cdot$\hss}}

\newcommand{\inorm}[1]{\mbox{$\left\Vert #1 \right\Vert_{\infty}$}}
\newcommand{\snorm}[1] {\Vert #1 \Vert}

\newcommand{\s}{\sigma}

\newtheorem{thm}{Theorem}

\theoremstyle{definition}

\newtheorem{rmk}[thm]{Remark}
\newtheorem{eg}{Example}
\newtheorem{pb}{Problem}

\hypersetup{
    linkcolor=blue,
}

\title{Towards Machine Wald}

\author{Houman Owhadi and Clint Scovel\\ California Institute of Technology}

\date{\today}

\renewcommand{\thefigure}{\arabic{figure}}

\makeatletter
\renewcommand{\p@subfigure}{\thefigure}
\makeatother

\newcounter{mycount}

\makeatletter
\def\blfootnote{\gdef\@thefnmark{}\@footnotetext}
\makeatother

\begin{document}

\maketitle

\begin{abstract}
The past century  has seen a steady increase in the need of estimating and predicting complex systems and making (possibly critical) decisions with limited information. Although computers have made possible the numerical evaluation of sophisticated statistical  models, these  models are still designed \emph{by humans} because there is currently no known recipe or algorithm for dividing the design of a statistical model into a sequence of arithmetic operations.
Indeed  enabling computers to \emph{think} as \emph{humans} have the ability to do when faced with uncertainty  is challenging in several major ways:
 (1) Finding optimal statistical models remains to be formulated as a well posed problem  when information on the system of interest is incomplete and comes in the form of  a complex combination of sample data, partial knowledge of constitutive relations and a limited description of the distribution of input random variables.
 (2) The space of admissible scenarios along with the space of relevant information, assumptions, and/or beliefs, tends to be infinite dimensional, whereas calculus on a computer is necessarily discrete and finite.
With this purpose, this paper explores the foundations of a rigorous framework for the scientific computation of optimal statistical estimators/models and reviews their connections with  Decision Theory, Machine Learning,
 Bayesian Inference, Stochastic Optimization, Robust Optimization, Optimal Uncertainty Quantification and Information Based Complexity.
\end{abstract}

%\blfootnote{\noindent ${^1}$Department of  Computing + Mathematical Sciences, California Institute of Technology, Pasadena CA 91125, USA. E-mail: owhadi@caltech.edu}
%\blfootnote{\noindent ${^2}$Department of  Computing + Mathematical Sciences, California Institute of Technology, Pasadena CA 91125, USA. E-mail: clintscovel@gmail.com}
\section{Introduction}
During the past century the need to solve large complex problems in applications such as fluid dynamics, neutron transport or ballistic prediction drove the parallel development of computers and numerical methods for solving ODEs and PDEs. It is now clear that this development lead to a paradigm shift. Before: each new PDE required the development of new theoretical methods and the employment of large teams of mathematicians and physicists; in most cases, information on solutions was only qualitative and based on general analytical bounds on fundamental solutions. After: mathematical analysis and computer science worked in synergy to give birth to robust numerical methods (such as finite element methods) capable of solving a large spectrum of PDEs without requiring the level of expertise of an
 A.~L.~Cauchy or level of insight of a R.~P.~Feynman. This transformation can be traced back to sophisticated calculations  performed by arrays of \emph{human computers} organized as parallel clusters such as in the pioneering work of  Lewis Fry Richardson \cite{Richardson1922, Lynch2008}, who
in 1922 had a room full of clerks attempt to solve finite-difference equations for the purposes of weather forecasting, and the 1947 paper by John Von Neumann and Herman Goldstine on Numerical Inverting of Matrices of High Order \cite{NeumannGoldstine1947}.
 Although Richardson's predictions failed due to the use of unfiltered data/initial conditions/equations and large time-steps not satisfying the CFL stability condition \cite{Lynch2008},  his vision was shared by Von Neumann
 \cite{Lynch2008}  in his proposal of the Meteorology Research Project to the U.S. Navy in 1946, qualified by Platzman \cite{Platzman1979} as ``perhaps the most visionary prospectus for numerical weather prediction since the publication of Richardson’s book a quarter-century earlier.''

The past century  has also seen a steady increase in the need of estimating and predicting complex systems and making (possibly critical) decisions with limited information. Although computers have made possible the numerical evaluation of sophisticated statistical models, these  models are still designed \emph{by humans} through the employment of multi-disciplinary teams of physicists, computer scientists and statisticians.
Contrary to the original \emph{human computers} (such as the ones pioneered by L. F. Richardson or overseen by R. P. Feynman at Los Alamos), these \emph{human teams} do not follow a specific algorithm (such as the one envisioned in Richardson's Forecast Factory where 64,000 human computers would have been working in parallel and at high speed to compute world weather charts \cite{Lynch2008}) because there is currently no known recipe or algorithm for dividing the design of a statistical model  into a sequence of arithmetic operations. Furthermore, while \emph{human computers} were given a specific PDE or ODE to solve, these \emph{human teams} are not given a well posed problem with a well defined notion of solution. As a consequence different \emph{human teams} come up with different \emph{solutions} to the design of the statistical model  along with
different estimates on uncertainties.

Indeed  enabling computers to \emph{think} as \emph{humans} have the ability to do when faced with uncertainty  is challenging in several major ways:  (1) There is currently no known recipe or algorithm for dividing the design of a statistical model into a sequence of arithmetic operations (2) Formulating the search for an optimal statistical estimator/model as a well posed problem is not obvious
 when information on the \emph{system of interest} is incomplete and comes in the form of
 a complex combination of sample data, partial knowledge of constitutive relations and a limited description of the distribution of input random variables.
 (3) The space of admissible scenarios along with the space of relevant information, assumptions, and/or beliefs, tend
to be infinite dimensional, whereas calculus on a computer is necessarily discrete and finite.

The purpose of this paper is to explore the foundations of a rigorous/rational framework for the scientific computation of optimal statistical estimators/models for complex systems and review their connections with  Decision Theory,
Machine Learning,  Bayesian Inference, Stochastic Optimization, Robust Optimization, Optimal Uncertainty Quantification and Information Based Complexity,  the most fundamental of these being the simultaneous emphasis on {\em computation}
and {\em performance} as in Machine Learning  initiated by Valiant \cite{Valiant}.

\section{The UQ problem without sample data}\label{subsecwithousd}

\subsection{\u{C}eby\u{s}ev, Markov and Kre\u{\i}n}\label{refsubsec11}
Let us start with a simple warm-up problem
\begin{pb}\label{pb:0}
Let $\mathcal{A}$ be the set of measures of probability on $[0,1]$ having mean less than $m \in (0,1)$.
Let $\mu^\dagger$ be an unknown element of $\mathcal{A}$ and let $a\in (m,1)$.  What is $\mu^\dagger[X\geq a]$?
\end{pb}

Observe that given the limited information on  $\mu^\dagger$, $\mu^\dagger[X\geq a]$ could a priori be any number in the interval $\big[\mathcal{L}(\mathcal{A}),\mathcal{U}(\mathcal{A})\big]$ obtained by computing the  sup (inf) of $\mu[X\geq a]$ with respect to  all possible candidates for $\mu^\dagger$, i.e.
\begin{equation}\label{eq:wup1}
\mathcal{U}(\mathcal{A}):=\sup_{\mu \in \mathcal{A}} \mu[X \geq a]
\end{equation}
and
\begin{equation*}\label{eq:wup1inf}
\mathcal{L}(\mathcal{A}):=\inf_{\mu \in \mathcal{A}} \mu[X \geq a]
\end{equation*}
where
\begin{equation*}
\mathcal{A}:=\big\{\mu \in \mathcal{M}([0,1])\mid \E_\mu[X]\leq m\big\}
\end{equation*}
and $\mathcal{M}([0,1])$ is the set of Borel probability measures on $[0,1]$.
It is easy to observe that the extremum of \eqref{eq:wup1} can be achieved only when $\mu$ is the weighted sum of a Dirac mass at $0$ and a Dirac mass at $a$. It follows that, although \eqref{eq:wup1} is an infinite dimensional optimization problem, it can be reduced to the simple one-dimensional optimization problem obtained by
letting  $p \in [0,1]$ denote the weight of the Dirac mass at $1$ and
$1-p$ the weight of the Dirac mass at $0$:  \emph{Maximize $p$ subject to $a p=m$}, producing the Markov bound $\frac{m}{a}$ as solution.

Problems such as \eqref{pb:0} can be traced back to \u{C}eby\u{s}ev
  \cite[Pg.~4]{Krein} ``Given: length, weight, position of the centroid and moment of inertia of a material rod with a density
varying from point to point. It is required to find the most accurate limits for the weight of a certain
segment  of this rod.'' According to  Kre\u{\i}n \cite{Krein}, although \u{C}eby\u{s}ev did solve this problem, it was his student
 Markov who supplied the proof in his thesis. See Kre\u{\i}n \cite{Krein}  for an account of the
 history of this subject  along with substantial contributions by Kre\u{\i}n.

\subsection{Optimal Uncertainty Quantification}

The generalization of the process described in Subsection \ref{refsubsec11} to complex systems involving imperfectly known functions and measures is
  the point of view of Optimal Uncertainty Quantification (OUQ)  \cite{OSSMO:2011, McKernsStrand:2011, KidaneI:2012, AdamsII:2012, Sullivan:2013, Kamgaouq:2014}.
Instead of developing sophisticated mathematical solutions, the OUQ approach is to develop optimization problems and reductions,
 so that their solution may be implemented on a computer, as in Bertsimas and Popescu's \cite{BertsimasPopescu:2005} convex optimization approach to
\u{C}eby\u{s}ev  inequalities, and the Decision Analysis framework of Smith \cite{Smith:1995}.

To present this generalization, for a topological space $\mathcal{X}$, let $\mathcal{F}(\mathcal{X})$ be the space of real-valued measurable functions and
$\mathcal{M}(\mathcal{X})$ be the set of Borel probability measures on $\mathcal{X}$.  Let $\mathcal{A}$ be an arbitrary subset of $\mathcal{F}(\mathcal{X}) \times \mathcal{M}(\mathcal{X})$, and let $\Phi \colon \mathcal{A} \to \R$ be a function producing a quantity of interest.

\begin{pb}\label{pbuqnodata}
Let $(f^\dagger,\mu^\dagger)$ be an unknown element of $\mathcal{A}$. What is $\Phi(f^\dagger,\mu^\dagger)$?
\end{pb}
 Therefore, in absence of sample data,
 in the context of this generalization one is interested in estimating $\Phi(f^\dagger,\mu^\dagger)$, where $(f^\dagger,\mu^\dagger)\in \mathcal{F}(\mathcal{X}) \times \mathcal{M}(\mathcal{X})$ corresponds to an \emph{unknown reality}:  the function $f^\dagger$ represents a \emph{response function} of interest, and $\mu^\dagger$ represents the probability distribution of the inputs of $f^\dagger$.  If $\mathcal{A}$ represents all that is known about $(f^\dagger,\mu^\dagger)$ (in the sense that $(f^\dagger,\mu^\dagger)\in \mathcal{A}$ and that any $(f, \mu)\in \mathcal{A}$ could, a priori, be $(f^\dagger,\mu^\dagger)$ given the available information) then \cite{OSSMO:2011} shows that the quantities
\begin{eqnarray}
	\label{eq:defma1}
	\mathcal{U}(\mathcal{A})&:=& \sup_{(f, \mu)\in \mathcal{A}} \Phi(f, \mu)\\
	\label{eq:defma2}
	\mathcal{L}(\mathcal{A})&:=& \inf_{(f, \mu)\in \mathcal{A}} \Phi(f, \mu)
\end{eqnarray}
determine the inequality
\begin{equation}
	\label{ineq_OUQ}
	\mathcal{L}(\mathcal{A}) \leq \Phi(f^\dagger,\mu^\dagger) \leq \mathcal{U}(\mathcal{A}),
\end{equation}
to be optimal given the available information $(f^\dagger,\mu^\dagger) \in \mathcal{A}$ as follows:  It is simple to see that the inequality \eqref{ineq_OUQ}  follows from the assumption that $(f^\dagger,\mu^\dagger) \in \mathcal{A}$.  Moreover, for any $\varepsilon >0$  there exists a $(f,\mu) \in \mathcal{A}$ such that
\[
	\mathcal{U}(\mathcal{A})-\varepsilon < \Phi(f,\mu) \leq \mathcal{U}(\mathcal{A}) .
\]
Consequently since all that is known about $(f^\dagger,\mu^\dagger)$ is that $(f^\dagger,\mu^\dagger) \in \mathcal{A}$, it follows that the upper bound $\Phi(f^\dagger,\mu^\dagger) \leq \mathcal{U}(\mathcal{A})$ is the best obtainable given that information, and the lower bound is optimal in the same sense.

Although the OUQ optimization problems \eqref{eq:defma1} and \eqref{eq:defma2} are extremely large and although some are computationally intractable,   an important subclass enjoys significant and practical finite-dimensional reduction properties \cite{OSSMO:2011}.  First, by \cite[Cor.~4.4]{OSSMO:2011}, although the optimization variables $(f, \mu)$ lie in a product space of functions and probability measures, for OUQ problems governed by linear inequality constraints on generalized moments, the search can be reduced to one over probability measures that are products of finite convex combinations of Dirac masses with explicit upper bounds on the number of Dirac masses.

Furthermore, in the special case that all constraints are generalized moments of functions of $f$, the dependency on the coordinate positions of the Dirac masses is eliminated by observing that the search over admissible functions reduces to a search over functions on an $m$-fold product of finite discrete spaces, and the search over $m$-fold products of finite convex combinations of Dirac masses reduces to a search over the products of probability measures on this $m$-fold product of finite discrete spaces \cite[Thm.~4.7]{OSSMO:2011}.  Finally, by \cite[Thm.~4.9]{OSSMO:2011}, using the lattice structure of the space of functions, the search over these functions can be reduced to a search over a finite set.

Fundamental to this development is Winkler's \cite{Winkler:1988}
 generalization of the characterization of the extreme points of compact (in the weak topology) sets of probability
measures constrained by a finite number of generalized moment inequalities defined by  continuous functions
  to {\em non-compact sets of tight measures},
in particular probability measures on Borel subsets of Polish metric spaces, defined by Borel measurable moment functions, along with his
\cite{winkler1978integral} development of Choquet theory for weakly closed convex {\em non-compact}
 sets of tight measures.  These results are based on  Kendall's \cite{Kendall}
 equivalence between  a linearly compact Choquet simplex and a
vector lattice  and
results of  Dubins \cite{Dubins} concerning the extreme points of affinely constrained
convex sets in terms of the extreme points of the unconstrained set. It is interesting to
note that Winkler \cite{Winkler:1988} uses Kendall's result to derive a strong sharpening of Dubins result \cite{Dubins}.
 Winkler's results allow the extension of existing optimization results over measures
on  compact metric spaces  constrained by continuous generalized moment functions
 to optimization over measures on Borel subsets of Polish spaces constrained by Borel measurable moment functions.
For systems with symmetry,  the Choquet theorem
of Varadarajan \cite{Varadarajan_groups} can be used to show that the Dirac masses can be replaced by
the ergodic measures in these results.
 The inclusion of sets of functions along with sets of measures in the optimization
problems facilitates  the application to systems with imprecisely known response functions. In particular,
  a result of Ressel \cite{Ressel_some}, providing conditions under which the map $(f,\mu) \rightarrow f_{*}\mu$
from function/measure pairs to the induced law is  Borel measurable,  facilitates  the extension of these techniques    from sets of {\em measures} to sets of {\em random variables}.  In general, the inclusion of functions
in the domain of optimization requires
 the development
of generalized  programming techniques such as generalized Benders decompositions described in
Geoffrion \cite{Geoffrion:Benders}.  Moreover, as has been so successful in Machine Learning,
it will be convenient to approximate the space of measurable functions $\mathcal{F}(\mathcal{X})$ by
some reproducing kernel Hilbert space $\mathcal{H}(\mathcal{X})\subset \mathcal{F}(\mathcal{X})$  producing
an approximation $\mathcal{H}(\mathcal{X})\times \mathcal{M}(\mathcal{X})\subset \mathcal{F}(\mathcal{X})
\times \mathcal{M}(\mathcal{X})$ to
the full base space. Under the mild assumption that $\mathcal{X}$ is an analytic subset of a Polish space and
 $\mathcal{H}(\mathcal{X})$ possesses a measurable
feature map, it has recently been shown in
\cite{OwhadiScovelSeparable}
that $\mathcal{H}(\mathcal{X})$ is separable. Consequently, since
all separable Hilbert spaces are isomorphic with $\ell^{2}$, it follows  that
the space $\ell^{2}\times \mathcal{M}(\mathcal{X})$ is a universal representation space
for the  approximation of $\mathcal{F}(\mathcal{X})
\times \mathcal{M}(\mathcal{X})$.  Moreover, in that case, since $\mathcal{X}$ separable and metric so is
$\mathcal{M}(\mathcal{X})$ in the weak topology, and since $\mathcal{H}(\mathcal{X})$ is Polish, it follows that
  the approximation
space $\mathcal{H}(\mathcal{X})\times \mathcal{M}(\mathcal{X})$ is the product of a Polish space
and a separable metric space. When furthermore $\mathcal{X}$ is Polish it follows
that the approximation space is the product of Polish spaces and therefore Polish.

\begin{eg}
	\label{eg:1}
	A classic example is $\Phi(f,\mu):=\mu[f\geq a]$ where $a$ is a safety margin. In the certification context one is interested in showing that $\mu^\dagger[f^\dagger\geq a]\leq \varepsilon$, where $\varepsilon$ is a safety certification threshold (i.e.\ the maximum acceptable $\mu^\dagger$-probability of the system $f^\dagger$ exceeding the safety margin $a$).  If $\mathcal{U}(\mathcal{A}) \leq \varepsilon$, then the system associated with $(f^\dagger,\mu^\dagger)$ is safe even in the worst case scenario (given the information represented by $\mathcal{A}$).  If $\mathcal{L}(\mathcal{A}) > \varepsilon$, then the system associated with $(f^\dagger,\mu^\dagger)$ is unsafe even in the best case scenario (given the information represented by $\mathcal{A}$).  If $\mathcal{L}(\mathcal{A}) \leq \varepsilon < \mathcal{U}(\mathcal{A})$, then the safety of the system cannot be decided (although one could declare the system to be unsafe due to lack of information).
\end{eg}

\subsection{Worst case analysis} The proposed solutions to Problems \ref{pb:0} and \ref{pbuqnodata} are particular instances of worst case analysis that,  as noted by \cite{Sniedovich2011} and \cite[p.~5]{RustemHowe:2002} is an old concept that could be summarized  by the popular adage \emph{When in doubt, assume the worst!} or
\begin{quote}
The gods to-day stand friendly, that we may,\\
Lovers of peace, lead on our days to age\!\\
But, since the affairs of men rests still uncertain,\\
Let’s reason with the worst that may befall.\\
\newline
Julius Caesar, Act 5, Scene 1\\
William Shakespeare (1564–-1616)
\end{quote}
As noted in \cite{OSSMO:2011}, an example of worst case analysis in  seismic engineering is that of  Drenick's \emph{critical excitation} \cite{Drenick:1973} which seeks to quantify the safety of a structure to the worst earthquake given a constraint on its magnitude. The combination of structural optimization (in various fields of engineering) to produce an optimal design given the (deterministic) worst-case scenario has been referred to as \emph{Optimization and Anti-Optimization} \cite{ElishakofOhsaki:2010}. The main difference between OUQ and Anti-optimization lies in the fact that the former is based on an optimization over (admissible) functions and measures $(f,\mu)$, while the latter only involves an optimization over $f$.
Because of its robustness, many engineers have adopted the (deterministic) worst-case scenario approach to UQ \cite[Chap.~10]{ElishakofOhsaki:2010} when a high  reliability is required.

\subsection{Stochastic and Robust Optimization}
  Robust Control \cite{Zhou96} and Robust Optimization \cite{Bentalbook2009, Bertsimas2011} have been founded upon the worst case approach to uncertainty. Recall that Robust Optimization describes optimization involving uncertain parameters. While these uncertain parameters are modeled as random variables (of known distribution) in Stochastic Programming \cite{Dantzig1955}, Robust Optimization only assumes that they are contained in known (ambiguity) sets.
Although, as noted in \cite{ElishakofOhsaki:2010}, \emph{probabilistic methods do not find appreciation among theoreticians and practitioners alike} because ``probabilistic reliability studies involve assumptions on the probability densities, whose knowledge regarding relevant input quantities is central'', the deterministic worst case approach (limited to optimization problems over $f$)  is sometimes ``too pessimistic to be practical'' \cite{Drenick:1973, ElishakofOhsaki:2010} because ``it does not take into account the \emph{improbability} that  (possibly independent or weakly correlated) random variables conspire to produce a failure event'' \cite{OSSMO:2011}  (which constitutes one motivation for considering ambiguity sets involving both measures and functions). Therefore OUQ and \emph{Distributionally Robust Optimization} (DRO) \cite{Bentalbook2009, Goh2010, Bertsimas2011, XuYu2014, Kuhn2013, Kuhn2014, hanasusanto2015distributionally} could be seen as middle-ground between the deterministic worst case approach of Robust Optimization \cite{Bentalbook2009, Bertsimas2011} and  approaches of  \emph{Stochastic Programming} and \emph{Chanced Constrained Optimization} \cite{Bot2010, Chen2010} by defining optimization objectives and constraints in terms of expected values and probabilities with respect to imperfectly known distributions.

Although  Stochastic Optimization has different objectives than OUQ and DRO,
many of its optimization results, such as found in Birge and Wets \cite{BirgeWets}, Ermoliev \cite{Ermoliev} and Gaivoronski, \cite{Gaivoronski}, are  useful. In particular,  the well-developed subject of
 Edmundson and Madansky  bounds such as
Edmundson \cite{Edmundson}, Madansky \cite{Madansky, Madansky_inequalities}, Gassman and Ziemba \cite{GassmannZiemba}, Huang, Ziemba and Ben-Tal  \cite{Huang}, Frauendorfer \cite{Frauendorfer}, Ben-Tal and Hochman \cite{Ben-Tal}, Huang, Vertinsky and Ziemba \cite{Huang_sharp},  and Kall \cite{Kall} provide powerful results.
 Recently Hanasusanto, Roitch, Kuhn  and Wiesemann \cite{hanasusanto2015distributionally}
 derive explicit conic reformulations for tractable problem classes and
suggest efficiently computable conservative approximations for intractable ones. In some cases, e.g.~Bertsimas and Popescu's \cite{BertsimasPopescu:2005} and   Han et al.~\cite{han2013convex},  DRO/OUQ optimization problems can  be  reduced to convex optimization.

\subsection{\u{C}eby\u{s}ev inequalities and optimization theory}
As noted in \cite{OSSMO:2011}, inequalities \eqref{ineq_OUQ} can be seen as a generalization of \u{C}eby\u{s}ev inequalities. The history of classical inequalities can be found in \cite{KarlinStudden:1966}, and some generalizations in \cite{BertsimasPopescu:2005} and \cite{VandenbergheBoydComanor:2007};  in the latter works, the connection between \u{C}eby\u{s}ev inequalities and optimization theory is developed  based on the work of Mulholland and Rogers \cite{MulhollandRogers:1958}, Godwin \cite{Godwin:1955}, Isii \cite{Isii:1959, Isii:1960, Isii:1963}, Olkin and Pratt \cite{OlkinPratt:1958}, Marshall and Olkin \cite{MarshallOlkin:1960}, and the classical Markov--Krein Theorem \cite[pages 82 \&~157]{KarlinStudden:1966}, among others.
 We also refer to the field of majorization, as discussed in Marshall and Olkin \cite{MarshallOlkin:1979}, the inequalities of Anderson \cite{Anderson:1955}, Hoeffding \cite{Hoeffding:1956b}, Joe \cite{Joe:1987}, Bentkus \emph{et al.}\ \cite{BentkusGeuzeZuijlen:2006}, Bentkus \cite{Bentkus:2002, Bentkus:2004}, Pinelis \cite{Pinelis:2007, Pinelis:2008}, and Boucheron, Lugosi and Massart \cite{BoucheronLugosiMassart:2000}. Moreover, the solution of the resulting nonconvex optimization problems benefit from duality theories for nonconvex optimization problems such as Rockafellar \cite{Rockafellar:1974} and the development of convex envelopes for them, as can be found, for example, in Rikun \cite{Rikun:1997} and Sherali \cite{Sherali:1997}.

\section{The UQ problem with sample data}

\subsection{From Game Theory to Decision Theory}
To motivate the general formulation in  the presence of sample data, consider another simple warm-up problem.
\begin{pb}\label{pb:1}
Let $\mathcal{A}$ be the set of measures of probability on $[0,1]$ having mean less than $m\in (0,1)$.
Let $\mu^\dagger$ be an unknown element of $\mathcal{A}$ and let $a\in (m,1)$. You observe $\Ddata:=(\Ddata_1,\ldots,\Ddata_n)$, $n$ i.i.d. samples from $\mu^\dagger$. What is the {\em sharpest} estimate of $\mu^\dagger[X\geq a]$?
\end{pb}
The only difference between problems \ref{pb:1} and \ref{pb:0} lies in the availability of data sampled from the underlying unknown distribution.
Observe that, in presence of this sample data, the notions of  \emph{sharpest estimate} or \emph{smallest interval of confidence}  are far from being transparent and call for  clear and precise definitions.
Note also that if the constraint $\E_{\mu^\dagger}[X]\leq m$ is ignored, and the number $n$ of sample data is large, then one could use the Central Limit Theorem or a concentration inequality (such as Hoeffding's inequality) to derive an interval of confidence for $\mu^\dagger[X\geq a]$. A non trivial question of practical importance is what to do when $n$ is not large.

Writing $\Phi(\mu^\dagger):=\mu^\dagger[X\geq a]$ as the quantity of interest, observe that an estimation of $\Phi(\mu^\dagger)$ is a function (which  will be written $\theta$) of the data $\Ddata$. Ideally one would like to choose $\theta$ so that the estimation error $\theta(\Ddata)-\Phi(\mu^\dagger)$ is as close as possible to zero. Since $\Ddata$ is random, a more robust notion of error is that of a statistical error $\Er\big(\theta, \mu^\dagger)$ defined by weighting
the error with respect to  a real  measurable positive loss function $V\colon\R\rightarrow \R$ and the distribution of the data, i.e.,
\begin{equation}\label{eq:stater0}
\Er(\theta, \mu^\dagger):=\E_{\Ddata \sim (\mu^\dagger)^n}\Big[V\big[\theta(\Ddata)-\Phi(\mu^\dagger)\big]\Big]
\end{equation}
Note that for $V(x)=x^2$, the statistical error $\Er(\theta, \mu^\dagger)$ defined in \eqref{eq:stater0} is the mean squared error with respect to the distribution of the data $\Ddata$ of the estimation error.
 For  $V(x)=\one_{[\gamma,\infty]}(|x|)$ defined for some $\gamma>0$, $\Er(\theta, \mu^\dagger)$ is the probability with respect to the distribution of $\Ddata$ that the estimation error is larger than $\gamma$.

Now since $\mu^\dagger$  is unknown, the statistical error $\Er(\theta, \mu^\dagger)$ of  any $\theta$
is also unknown. However one can still identify the least upper bound on that statistical error through a worst case scenario with respect to all possible candidates for $\mu^\dagger$, i.e.
\begin{equation}\label{eq:stater1}
\sup_{\mu \in \mathcal{A}}\Er(\theta, \mu)\, .
\end{equation}
The sharpest estimator (possibly within a given class) is then naturally obtained as the minimizer of \eqref{eq:stater1} over all functions $\theta$ of the data $\Ddata$ within that class, i.e. as the minimizer of
\begin{equation}\label{eq:stater2}
\inf_{\theta}\sup_{\mu \in \mathcal{A}}\Er\big(\theta, \mu)\, .
\end{equation}

 Observe that the optimal estimator is identified independently from  the observation/realization of  the data and if the minimum of \eqref{eq:stater2} is not achieved then one can still use a near-optimal $\theta$. Then, when the data is observed, the   estimate of the quantity of interest $\Phi(\mu^\dagger)$ is then derived by evaluating the near-optimal $\theta$ on the data $\Ddata$.

The notion of optimality described here  is that of Wald's Statistical Decision Theory  \cite{Wald:1939, Wald:1945, Wald1947, Wald49, Wald:1950}, evidently influenced by Von Neumann's Game Theory \cite{VNeumann28, VonNeumann:1944}. In Wald's formulation \cite{Wald:1945}, which cites both Von Neumann \cite{VNeumann28}
and Von Neumann and Morgenstern \cite{VonNeumann:1944}, the statistician finds himself in an adversarial game played against the Universe in which he tries to minimize a risk function $\Er\big(\theta, \mu)$ with respect to $\theta$ in a worst case scenario with respect to what the Universe's choice of $\mu$ could be.

\subsection{The optimization approach to statistics} The optimization approach to statistics is not new and
this section will now give a short, albeit incomplete, description of its development, primarily using Lehmann's account
   \cite{lehmann2009some}. Accordingly, it began with
Gauss and Laplace with the non-parametric result referred to as the Gauss-Markov theorem, asserting that
the least squares estimates are the linear unbiased estimates with minimum variance.
Then, in    Fisher's fundamental paper \cite{fisher1922mathematical}, for parametric models he proposes the maximum
likelihood estimator and claims (but does not prove) that such estimators are consistent and asymptotically efficient.
According to Lehmann, ``the situation is complex, but under suitable restrictions Fisher's conjecture is essentially correct ..."  Fisher's maximum likelihood principle was first proposed on  intuitive grounds and then its optimality properties developed. However, according to Lehmann \cite[Pg.~1011]{lehmann:optimality_symposia},
 Pearson then asked Neyman ``Why these tests rather than any of the many other that could be proposed?
This question resulted in
Neyman and Pearson's  1928 papers \cite{neyman1928use} on the likelihood ratio method, which  gives the same answer as Fisher's tests under normality assumptions.  However,
 Neyman was not satisfied. He agreed that the likelihood ratio principle was appealing but
 felt that it was lacking a logically convincing justification. This then led  to the publication Neyman and Pearson \cite{neyman1933problem}, containing
their now famous Neyman-Pearson Lemma, which according to Lehmann \cite{lehmann2009some}, ``In a certain sense this is the true start of optimality theory."
In a major extension of the Neyman-Pearson work, Huber \cite{Huber:1964} proves a {\em robust} version of the Neyman-Pearson Lemma, in particular, providing an optimality criteria defining the robust estimator,  giving rise to a rigorous theory of robust statistics based on optimality, see Huber's Wald Lecture \cite{Huber:Wald}. Robustness is particularly suited to the Wald framework since robustness considerations are easily formulated with the proper choices of admissible functions and measures
in the Wald framework. Another example is Kiefer's introduction of optimality into experimental design,   resulting
in
 Kiefer's 40 papers on  Optimum Experimental Designs \cite{kiefer1985collected}.

 Not everyone was happy with ``optimality'' as a guiding principle. For example Lehmann \cite{lehmann2009some} states
that at a 1958 meeting of the Royal Statistical Society at which Kiefer presented a survey talk \cite{kiefer1959optimum}
on Optimum Experimental Designs,  Barnard quotes Kiefer
as saying of procedures proposed in a paper by Box and Wilson that they ``often [are] not even
well-defined rules of operation." Barnard's reply:
\begin{quotation}
        \noindent
 ``in the field of practical human activity, rules of operation which are not well-defined may be preferable to rules which are."
\end{quotation}
  Wynn \cite{wynn1992introduction}, in his  introduction to a reprint of Kiefer's paper,
 calls this ``a clash of statistical cultures."
Indeed, it is interesting to read the generally negative responses to  Kiefer's article \cite{kiefer1959optimum}
 and the
remarkable rebuttal by Kiefer therein.
 Tukey had other criticisms regarding ``The tyranny of the best" in \cite{tukey1961statistical} and
 ``The dangers of optimization"  in
\cite{tukey1962future}. In the latter he writes:
\begin{quotation}
        \noindent ``Some [statisticians] seem to equate [optimization] to statistics an attitude which,
if widely adopted, is guaranteed to produce a dried-up, encysted field with little chance of real growth."
\end{quotation}
For an account of how the Exploratory Data Analysis
approach of Tukey fits within the Fisher/Neyman-Pearson debate,
 see Lehnard \cite{lenhard2006models}.

Let us also remark on the influence that Student -William Sealy Gosset
 had on both Fisher and Pearson.
 As presented in Lehmann's \cite{Lehmann_Student} ````{S}tudent" and small-sample theory",
 quoting F.~N.~David \cite{laird1989conversation}:
``I think he [Gosset] was really the big influence
in statistics... He asked the questions and Pearson or Fisher put them into statistical language and then Neyman
came to work with the mathematics. But I think most of it stems from Gosset."
The aim of Lehmann's paper \cite{Lehmann_Student} is to consider to what extend David's conclusion is justified.
Indeed, the claim is surprising since Gosset is mainly known for only one contribution, that is Student
 \cite{student1908probable}, with the introduction of Student's t-test and its analysis under the normal distribution.
According to Lehmann ``Today the pathbreaking nature of this paper is generally recognized and has been widely commented upon, \ldots''.
Gosset's primary concern in communicating with both Fisher and Pearson  was the robustness of the test to non-normality. Lehmann concludes that ``the main ideas
leading to Pearson's research were indeed provided by Student."  See Lehmann \cite{Lehmann_Student}
for the full account, including Gosset's relationship to the
 Fisher-Pearson debate,
  Pearson \cite{Pearson} for a Statistical biography of Gosset, and Fisher \cite{fisher1939student} for a eulogy.  Consequently, modern statistics appears to owe
a lot to Gosset. Moreover, the reason for the pseudonym was a
policy by Gosset's employer, the brewery Arthur
Guinness Sons and Co., against work done for the
firm being made public. Allowing Gosset to publish
under a pseudonym was a concession that resulted
in the birth of the statistician Student.
 Consequently, the authors would like to take this opportunity to thank the Guinness Brewery  for
its influence on statistics today, and for such a fine beer.

\subsection{Abraham Wald}
Following Neyman and Pearson's breakthrough,
 a different approach to optimality was introduced in
  Wald \cite{Wald:1939} and then  developed in a sequence of papers   culminating
 in Wald's \cite{Wald:1950} book Statistical Decision Functions.
 Evidence of the influence of Neyman on Wald
can be found in the citation of  Neyman \cite{neyman1937outline} in  the introduction of Wald \cite{Wald:1939}.
Brown \cite{brown2000essay} quotes the students of Neyman in 1966 from \cite{neyman1967selection}:
\begin{quotation}
        \noindent
``The concepts of confidence intervals and of the Neyman-Pearson
theory have proved immensely fruitful. A natural
but far reaching extension of their scope can be found in
Abraham Wald's theory of statistical decision functions. The
elaboration and application of the statistical tools related to
these ideas has already occupied a generation of statisticians.
It continues to be the main lifestream of theoretical statistics."
\end{quotation}
Brown's purpose was to address  if the last sentence
in the quote  was still true in 2000.

 Wolfowitz \cite{wolfowitz1952abraham} describes the primary accomplishments of Wald's Statistical Decision Theory as follows:
\begin{quotation}
        \noindent
``Wald's greatest achievement was the theory of statistical decision functions, which includes almost all
problems which are the raison d'etre of statistics"
\end{quotation}
  Leonard \cite[Chp.~12]{leonard2010neumann} portrays  Von Neumann's return to Game Theory as ``partly an early reaction to
upheaval and war". However he adds that eventually Von Neumann became personally involved in the war effort and
``with that involvement came a significant, unforeseeable moment in the history of game theory: this new mathematics
made its wartime entrance into the world, not as the abstract theory of social order central to the book, but
as a problem solving technique."
Moreover, on pages
 278--280  Leonard discusses the   Statistical Research Groups at Berkeley, Columbia, and Princeton, in particular
Wald at Columbia, and how the effort to develop inspection and testing procedures leads Wald  to the development
of sequential methods
 ``apparently yielding significant economies in inspection in the Navy",
leading to the publication of Wald and Wolfowitz' \cite{Wald:optimum} proof of
the  optimality   of the sequential probability ratio test and Wald's book \cite{Wald:sequential} Sequential Analysis.
Leonard's claim essentially is that the war stimulated these fine theoretical minds to pursue
activities with real application value. In this regard, it is relevant to note
Mangel and Samaniego's \cite{Mangel} stimulating description of
Wald's work on aircraft survivability, along with the  contemporary, albeit somewhat vague,  description of
``How a Story from World War II shapes Facebook today" by
  Wilson \cite{Wilson}.   Indeed,
in the problem of how to
 allocate armoring to the allied bombers based on their condition upon return from their missions,
it was discovered that armoring where the previous planes had been hit was not improving their rate of return.
Wald's ingenious insight was that these were the returning bombers not the ones which had been shot down. So
the places where the returning bombers were hit are more likely to be the places where one {\em does not}
need to add armoring. Evidently, his rigorous and unconventional innovations to transform this intuition
into a real methodology saved many lives.
 Wolfowitz \cite{wolfowitz1952abraham} states
\begin{quotation}
        \noindent
``Wald not only posed his statistical problems clearly and
precisely, but he posed them to fit the practical problem and to accord with the
decisions the statistician was called on to make. This, in my opinion, was the
key to his success-a high level of mathematical talent of the most abstract
sort, and a true feeling for, and insight into, practical problems. The combination
of the two in his person at such high levels was what gave him his outstanding
character."
\end{quotation}

The section on Von Neumann and Remak (along with the Notes that follows it) in
Kurz and Salvadori \cite{kurz2002understanding}
describes  Wald and Von Neumann's relations.
Brown \cite{brown1994minimaxity}  credits Wald as the creator  of the
 minmax idea  in statistics,
evidently given axiomatic justification by Gilboa and Schmeidler \cite{GilboaSchmeidler_maxmin}.
This  certainly had something to do with his
 friendship with Morgenstern  and his relationship with Von Neumann, who together authored
 the famous book \cite{VonNeumann:1944}, but this influence can be explicitly seen
in Wald's \cite{Wald:1945} citation of Von Neumann \cite{VNeumann28}
and Von Neumann and Morgenstern \cite{VonNeumann:1944} in his
 introduction \cite{Wald:1945} of the minmax idea in Statistical Decision Theory.
 Indeed, Wolfowitz states that
\begin{quotation}
        \noindent
``...
he was also spurred on by the connection
between the newly announced results of [Von Neumann and Morgenstern] \cite{VonNeumann:1944} and his own theory, and by
the general interest among economists and others aroused by the theory of
games."
\end{quotation}
Wolfowitz asserts that Wald's work \cite{Wald:1939}
  Contributions to the Theory of Statistical Estimation and Testing Hypotheses
is ``probably his most important paper" but that it ``went almost completely unnoticed",  possibly
because ``The use of Bayes solutions was deterrent" and ``Wald did not really emphasize  that he was using Bayes
 solutions only as a tool."
Moreover, although Wolfowitz considered  Wald's Statistical Decision Functions \cite{Wald:1950} his
 greatest achievement, also says
\begin{quotation}
        \noindent
``The statistician who wants to apply the results of \cite{Wald:1950} to specific problems
is likely to be disappointed. Except for special problems, the complete classes
are difficult to characterize in a simple manner and have not yet been characterized.
Satisfactory general methods are not yet known for obtaining minimax
solutions. If one is not always going to use a minimax solution (to which.serious
objections have been raised) or a solution satisfying some given criterion, then
the statistician should have the opportunity to choose from among "representative"
decision functions on the basis of their risk functions. These are not
available except for the simplest cases. It is clear that much remains to be done
before the use of decision functions becomes common. The theory provides a
rational basis for attacking almost any statistical problem, and, when some
computational help is available and one makes some reasonable compromises
in the interest of computational feasibility, one can obtain a practical answer
to many problems which the classical theory is unable to answer or answers in
an unsatisfactory manner."
\end{quotation}

Wolfowitz \cite{wolfowitz1952abraham},  Morgenstern \cite{morgenstern1951abraham}
 and Hotelling \cite{hotelling1951abraham} provide a description of Wald's impact at the time of his passing. The influence of Wald's minimax paradigm can also be observed on (1) decision making under severe uncertainty \cite{Sniedovich2007, Sniedovich2011, Sniedovich2012} (2) Stochastic Programming \cite{Shapiro2002} (minimax analysis of stochastic problems) (3)  Minimax solutions of stochastic linear programming  problems \cite{Zavokva1966} (4) Robust Convex Optimization \cite{BenTalNemirovski1998} (where one must find the best decision in view of the worst case
parameter values within deterministic uncertainty sets) (4) Econometrics  \cite{Tintner1952} and (5) Savage's minimax Regret model \cite{Savage1951}.

\subsection{Generalization to unknown pairs of functions and measures and to arbitrary sample data.}
In practice, complex systems of interest may involve, both an imperfectly known response function $f^\dagger$ and an imperfectly known  probability measure $\mu^\dagger$ as illustrated in the following problem.
\begin{pb}\label{pb:3}
Let $\mathcal{A}$ be a set of real functions and measures of probability on $[0,1]$ such that $(f,\mu)\in \mathcal{A}$ if and only if $\E_{\mu}[X]\leq m$ and $\sup_{x\in [0,1]}\big|g(x)-f(x)\big|\leq 0.1$ where $g$ is some given real function on $[0,1]$. Let $(f^\dagger,\mu^\dagger)$ be an unknown element of $\mathcal{A}$ and let $a\in \R$. Let $(X_1,\ldots,X_n)$ be $n$ i.i.d. samples from $\mu^\dagger$, you observe $(d_1,\ldots,d_n)$ with $d_i=\big(X_i,f^\dagger(X_i)\big)$. What is the ``sharpest'' estimate of $\mu^\dagger\big[f(X)\geq a\big]$?
\end{pb}

Problem \ref{pb:3} is an illustration of a situation in which  the response function $f^\dagger$ and the probability measure $\mu^\dagger$ are not directly observed and the sample data arrives in the form of realizations of random variables, the distribution of which is related to $(f^\dagger,\mu^\dagger)$.  To simplify the current presentation  assume that this relation is, in general, determined by a function of $(f^\dagger,\mu^\dagger)$ and use the following notation:  $\Dspace$ will denote the observable space (i.e.\ the space in which the sample data $d$ take values, assumed to be a metrizable Suslin space) and $d$ will denote the $\Dspace$-valued random variable corresponding to the observed sample data.  To represent the dependence of the distribution of
$d$ on the unknown state $(f^\dagger,\mu^\dagger) \in \mathcal{A}$ introduce a measurable function
\begin{equation}\label{eqhhshduhd}
	\Dmap \colon \mathcal{A} \to \mathcal{M}(\Dspace),
\end{equation}
where $\mathcal{M}(\Dspace)$ is given the Borel structure corresponding to the weak topology, to define this relation.  The idea is that $\Dmap(f, \mu)$ is the probability distribution of the observed sample data $d$ if $(f^\dagger,\mu^\dagger)=(f, \mu)$, and for this reason it may be called the \emph{data map} (or, even more loosely, the \emph{observation operator}). Now consider the following problem.

\begin{pb}\label{pb:4}
Let $\mathcal{A}$ be a known subset of $\mathcal{F}(\mathcal{X}) \times \mathcal{M}(\mathcal{X})$ as in Problem \ref{pbuqnodata} and let $\Dmap$ be a known data map as in  \eqref{eqhhshduhd}.
Let $\Phi$ be a known measurable semi-bounded function mapping $\mathcal{A}$ onto $\R$.
Let $(f^\dagger,\mu^\dagger)$  be an unknown element of $\mathcal{A}$.  You observe $d\in \Dspace$  sampled from the distribution $\Dmap(f^\dagger,\mu^\dagger)$. What is the sharpest estimation of $\Phi(f^\dagger,\mu^\dagger)$?
\end{pb}

\subsection{Model error and optimal models}
As in Wald's Statistical Decision Theory  \cite{Wald:1945}, a natural notion of optimality can be obtained by formulating Problem \ref{pb:4} as an adversarial game in which player A chooses $(f^\dagger,\mu^\dagger)\in \mathcal{A}$ and player B (knowing $\mathcal{A}$ and $\Dmap$) chooses a function $\theta$ of the observed data $d$. As in \eqref{eq:stater0} this notion of optimality requires the introduction of a risk function
\begin{equation}\label{eqscoseloss}
\Er\big(\theta, (f,\mu)\big):=\E_{\Ddata \sim \Dmap(f,\mu)}\Big[V\big[\theta(\Ddata)-\Phi(f,\mu)\big]\Big]
\end{equation}
where  $V\colon\R\rightarrow \R$ is a real  positive measurable loss function function. As in \eqref{eq:stater1} the
 least upper bound on that statistical error $\Er\big(\theta, (f,\mu)\big)$ is obtained as  through a worst case scenario with respect to all possible candidates for $(f,\mu)$ (player's A choice), i.e.
\begin{equation}\label{eq:stater1scose}
\sup_{(f,\mu) \in \mathcal{A}}\Er\big(\theta, (f,\mu)\big)
\end{equation}
and an optimal estimator/model (possibly within a given class) is then naturally obtained as the minimizer of \eqref{eq:stater1scose} over all functions $\theta$ of the data $\Ddata$ in that class (player's B choice), i.e. as the minimizer of
\begin{equation}\label{eq:stater2scose}
\inf_{\theta}\sup_{(f,\mu) \in \mathcal{A}}\Er\big(\theta, (f,\mu)\big)\, .
\end{equation}
 Since in real applications true optimality will never be achieved,
 it is natural to generalize to considering near-minimizers of \eqref{eq:stater2scose} as near-optimal models/estimators.

\begin{rmk}
In situations where the data map is imperfectly known (e.g. when the data $d$ is corrupted by some noise of imperfectly known distribution) one has to include a supremum over all possible candidates $\Dmap\in \Dmaps$ in the calculation of the least upper bound on the statistical error.
\end{rmk}

\subsection{Mean squared error, variance and bias} For $(f,\mu)\in \mathcal{A}$ write  $\operatorname{Var}_{\Ddata \sim \Dmap(f,\mu)}\big[\theta(\Ddata)\big]$ the variance of the random variable $\theta(\Ddata)$ when $\Ddata$ distributed according to $\Dmap(f,\mu)$, i.e.
\begin{equation*}
\operatorname{Var}_{\Ddata \sim \Dmap(f,\mu)}\big[\theta(\Ddata)\big]:=\E_{\Ddata\sim \Dmap(f,\mu)}\Big[\big(\theta(\Ddata)\big)^2\Big]-\Big[\E_{\Ddata\sim \Dmap(f,\mu)}\big[\theta(\Ddata)\big]\Big]^2
\end{equation*}
The following equation, whose proof is straightforward, shows that for $V(x)=x^2$,
 the least upper bound on the mean squared error of $\theta$ is equal to the least upper bound on the sum of the variance of $\theta$ and the square of its bias.
\begin{equation*}
\sup_{(f,\mu)\in \mathcal{A}} \Er\big(\theta,(f,\mu)\big)= \sup_{(f,\mu)\in \mathcal{A}} \Bigg[ \operatorname{Var}_{\Ddata \sim \Dmap(f,\mu)}\big[\theta(\Ddata)\big]+
\Big(\E_{\Ddata\sim \Dmap(f,\mu)}\big[\theta(\Ddata)\big]-\Phi(f,\mu)\Big)^2 \Bigg]
\end{equation*}
Therefore, for $V(x)=x^2$, the bias/variance tradeoff is made explicit.

\subsection{Optimal interval of confidence}
Although $\Er$ can a priori be defined to be any  risk function, taking $V(x)=\one_{[\gamma,\infty]}(|x|)$ (for some $\gamma>0$) in \eqref{eq:stater0} allows for a transparent and objective identification of optimal intervals of confidence.
Indeed, writing,
\begin{equation*}\label{eq:stater0oi}
\Er_\gamma\big(\theta, (f,\mu)):=\P_{\Ddata \sim \Dmap(f,\mu)}\Big[\big|\theta(\Ddata)-\Phi(f,\mu)\big|\geq \gamma\Big]
\end{equation*}
note that $\sup_{(f,\mu) \in \mathcal{A}}\Er_\gamma\big(\theta, (f,\mu))$ is the least upper bound on the
probability (with respect to the distribution of $\Ddata$) that the difference between the true value of the quantity of interest $\Phi(f^\dagger,\mu^\dagger)$ and its estimated value  $\theta(\Ddata)$ is larger than $\gamma$. Let $\epsilon \in [0,1]$. Define
\begin{equation*}
\gamma_{\epsilon}:=\inf \big\{\gamma >0\mid  \inf_{\theta}\sup_{(f,\mu) \in \mathcal{A}}\Er_\gamma\big(\theta, (f,\mu)) \leq \epsilon \big\},
\end{equation*}
and observe that if $\theta_\epsilon$ is a minimizer of $\inf_{\theta}\sup_{(f,\mu) \in \mathcal{A}}\Er_{\gamma_\epsilon}\big(\theta, (f,\mu)\big)$ then
 $[\theta_\epsilon(\Ddata)-\gamma_{\epsilon}, \theta_\epsilon(\Ddata)+\gamma_{\epsilon} ]$ is the smallest interval of confidence (random interval obtained as a function of the data) containing $\Phi(f^\dagger,\mu^\dagger)$ with probability at least $1-\epsilon$.
Observe also that this formulation is a natural extension of the OUQ formulation as described in \cite{OSSMO:2011}. Indeed, in the absence of sample data, it is easy to show that $\theta_1$ is the midpoint of the optimal interval $[\mathcal{L}(\mathcal{A}),\mathcal{U}(\mathcal{A})]$.

\begin{rmk}
We refer to \cite{Fisher1935, Fisher1955, Spanos14} and in particular to Stein's notorious paradox \cite{Stein1956} for the importance of a careful choice for loss function.
\end{rmk}

\subsection{Ordering the space of experiments}
 A natural
  objective of UQ and Statistics is the design of experiments,
their comparisons and the identification of optimal ones. Introduced in Blackwell
 \cite{blackwell1953equivalent} and Kiefer
 \cite{kiefer1959optimum},  with a more modern perspective in Le Cam \cite{LeCam} and Strasser \cite{strasser1985mathematical}, here
observe that \eqref{eq:stater2scose}, as a function of $\Dmap$, induces an order (transitive, total, but not antisymmetric) on the space of data maps that has a natural experimental design interpretation. More precisely if the data maps  $\Dmap_1$ and $\Dmap_2$ are interpreted as the distribution of the outcome of two possible experiments and if the value of \eqref{eq:stater2scose} is smaller for $\Dmap_2$ than $\Dmap_1$, then $\Dmap_2$ is a preferable experiment.

\subsection{Mixing models} 	Given  estimators $\theta_1, \dots, \theta_n$ can one obtain a better estimator by mixing those estimators? If $V$ is convex (or quasi-convex) then	the problem of finding an $\alpha\in [0,1]^n$ minimizing the statistical error of $\sum_{i=1}^n \alpha_i \theta_i$ under the constraint $\sum_{i=1}^n \alpha_i=1$ is a finite-dimensional convex optimization problem in $\alpha$. If estimators are seen as models of reality then this observation supports the idea that one can obtain improved models by mixing them (with the goal of achieving minimal statistical errors).

\section{The Complete Class Theorem and  Bayesian inference}\label{sec_staterror_random}

\subsection{The Bayesian approach}
The Bayesian answer to Problem \ref{pb:4} is to assume that $(f^\dagger,\mu^\dagger)$  is a sample from some (prior) measure $\pi$ on $\mathcal{A}$ and then condition the expectation of $\Phi(f,\mu)$ with respect to the observation of the data, i.e. use
 \begin{equation}\label{eqcondexp}
\E_{(f,\mu)\sim \pi, d\sim \Dmap(f,\mu)}\big[\Phi(f,\mu)\big| d\big]
 \end{equation}
 as the estimator $\theta(d)$. This requires giving $\mathcal{A}$ the structure of a measurable space such that important quantities of interest such as $(f,\mu)\rightarrow \mu[f(X)\geq a]$ and $(f,\mu)\rightarrow \E_\mu[f]$ are measurable. This can be achieved using results of Ressel \cite{Ressel_some} providing conditions under which the map $(f,\mu) \rightarrow f_{*}\mu$
from function/measure pairs to the induced law is  Borel measurable. We will henceforth assume
$\mathcal{A}$ to be a Suslin space, and proceed to construct the measure of probability
 $\pi\odot \Dmap$ of $\big((f,\mu),d\big)$ on $\mathcal{A} \times \Dspace$
 via a natural generalization of the Campbell measure and Palm distribution
associated with a random measure as described in \cite{Kallenberg:1975}, see also \cite[Ch.~13]{Daley} for a more current treatment.  We refer to Subsection \ref{subsecpidmap} of the appendix for the details of the construction of
 the distribution $\pi \odot \Dmap$ of $\big((f,\mu),d\big)$ when
$(f,\mu)\sim \pi$ and $d\sim \Dmap(f,\mu)$, and of  the marginal distribution  $\pi \cdot \Dmap$
 of $\pi \odot \Dmap$ on the data space $\Dspace$,
 and the resulting regular conditional expectation \eqref{eqcondexp}. Consequently,
  the nested  expectation $\E_{(f,\mu)\sim \pi, d\sim \Dmap(f,\mu)}$ appearing in \eqref{eqcondexp}
 will from now on be rigorously written  as the expectation $\E_{((f,\mu),d)\sim \pi\odot \Dmap}$.

\paragraph{Statistical error when $(f^\dagger,\mu^\dagger)$ is random.}
When  $(f^\dagger,\mu^\dagger)$ is a random realization of $\pi^\dagger$ one may consider averaging the statistical error \eqref{eqscoseloss} with respect to $\pi^\dagger$ and introduce
\begin{equation}\label{eq:dsicee3rep0}
\Er(\theta, \pi^\dagger):=\E_{((f,\mu),\Ddata) \sim \pi^\dagger \odot \Dmap}\Big[V\big[\theta(\Ddata)-\Phi(f,\mu)\big]\Big]
\end{equation}

When $\pi^\dagger$ is an unknown element of a subset $\Pi$ of $\mathcal{M}(\mathcal{A})$ the least upper bound on the statistical error \eqref{eq:dsicee3rep0} given the available information is obtained by taking the sup of \eqref{eq:dsicee3rep0} with respect to all possible candidates for $\pi^\dagger$, i.e.
\begin{equation}\label{eq:dsicee3rep0pidmaps}
\sup_{\pi \in \Pi}\Er(\theta, \pi)
\end{equation}
When $\mathcal{A}$ is Suslin and
when $(f^\dagger,\mu^\dagger)$ is not a random sample from $\pi^\dagger$ but simply an unknown element of $\mathcal{A}$, then a straightforward application of the reduction theorems of \cite{OSSMO:2011} implies that
 when $\Pi=\mathcal{M}(\mathcal{A})$, that
 \eqref{eq:dsicee3rep0pidmaps}  is equal to \eqref{eq:stater2scose}, i.e.
\begin{equation}\label{eq:aberrswedesation0}
\sup_{(f,\mu)\in \mathcal{A}} \Er\big(\theta, (f,\mu)\big)= \sup_{\pi \in \mathcal{M}(\mathcal{A})}\Er(\theta, \pi)
\end{equation}

\subsection{Relation between adversarial model error and Bayesian error.}
When $\Phi$ has a second moment with respect to $\pi$, one can utilize  the classical variational description of conditional expectation as follows:
Denote by $L^2_{\pi \cdot \Dmap}(\Dspace)$ the space of ($\pi \cdot \Dmap$ a.e.~equivalence classes of) real-valued measurable functions on $\Dspace$ that are square-integrable with respect to the measure $\pi \cdot \Dmap$.
Then
 one has (see Subsection \ref{subseccondexp})
\begin{equation*}\label{eq:condl2pi}
	\E_{\pi \odot \Dmap}\big[ \Phi \big| d\big]:= \argmin_{h \in L^2_{\pi\cdot \Dmap}(\Dspace)}
	\E_{(f,\mu,d) \sim \pi \odot \Dmap}\Big[ \big(\Phi(f,\mu)-h(d)\big)^2 \Big]\,  .
\end{equation*}
In other words, if $(f,\mu)$ is sampled from the measure $\pi$, $\E_{\pi \odot \Dmap}\big[ \Phi(f,\mu)\big| d\big]$ is the best mean-square approximation of $\Phi(f,\mu)$ in the space of square-integrable functions of $d$.  As with the regular conditional probabilities, the  real-valued function on $\Dspace$
\begin{equation}\label{eq:dsicees3repu}
\theta_{\pi}(\Ddata)=\E_{(f,\mu,\Drv)\sim\pi \odot \Dmap}\big[\Phi(f,\mu)\big| \Drv=\Ddata\big],\quad \Ddata
\in \Dspace
\end{equation}
is uniquely defined up to subsets of $\Dspace$ of $(\pi\cdot \Dmap)$-measure zero.

Using the orthogonal projection property of the conditional expectation one obtains that if $V(x)=x^2$, then for arbitrary $\theta$,
\begin{equation}\label{eq:minawswaeerration0}
\begin{split}
\Er (\theta,\pi)=\Er(\theta_{\pi},\pi)+
\E_{\Ddata\sim \pi \cdot \Dmap}\Big[\theta(\Ddata)-\theta_{\pi}(\Ddata)\Big]^2
\end{split}
\end{equation}
Therefore, if $\Pi \subset \mathcal{M}(\mathcal{A})$   is an admissible set of priors, then \eqref{eq:minawswaeerration0} implies that
\begin{equation*}\label{eqminaswsuiswueeerration0a}
\inf_{\theta} \sup_{\pi\in \Pi}\Er(\theta,\pi)\geq \sup_{\pi\in \Pi} \Er(\theta_{\pi},\pi)\, .
\end{equation*}
In particular,  when $\Pi = \mathcal{M}(\mathcal{A})$ \eqref{eq:aberrswedesation0} implies
that
\begin{equation}\label{eqminaswsuiswueeerration0}
\inf_{\theta} \sup_{(f,\mu)\in \mathcal{A}}\Er(\theta,(f,\mu))\geq \sup_{\pi\in \mathcal{M}(\mathcal{A})} \Er(\theta_{\pi},\pi)\, .
\end{equation}
Therefore, the mean squared error of the best estimator assuming $(f^\dagger,\mu^\dagger) \in \mathcal{A}$
 to be unknown
  is bounded  below by the largest  mean squared error
 of the Bayesian estimator  obtained by assuming that $(f^\dagger,\mu^\dagger)$  is distributed
according to some $\pi \in \mathcal{M}(\mathcal{A})$.
In the next section it will be shown that
  a complete class theorem can be used to obtain that  \eqref{eqminaswsuiswueeerration0} is actually an equality.
In  that case,
 \eqref{eqminaswsuiswueeerration0} can be used to quantify the approximate optimality of an estimator  by comparing the least upper bound $\sup_{(f,\mu)\in \mathcal{M}(\mathcal{A})} \Er\big(\theta, (f,\mu)\big)$ on the error of that estimator  with $\Er(\theta_{\pi},\pi)$ for a carefully chosen $\pi$.

\subsection{Complete Class theorem}
A fundamental question is whether \eqref{eqminaswsuiswueeerration0} is an equality: is the adversarial  error of the best estimator equal to the non-adversarial error of the worst Bayesian estimator? Is the best estimator Bayesian or an approximation thereof?
A remarkable result of decision theory \cite{Wald:1939, Wald:1945, Wald1947, Wald49, Wald:1950} is the complete class theorem which states (in the formulation of this paper) that if (1) the admissible measures $\mu$ are absolutely continuous with respect the Lebesgue measure (2)
 the
loss function $V$ in the definition of $\Er\big(\theta,(f,\mu)\big)$
 is convex in $\theta$ and (3) the decision space is compact, then optimal estimators live in the Bayesian class and non-Bayesian estimators cannot be optimal. The idea of the proof of this result is to use the compactness of the decision space and the continuity of the loss function to approximate the  decision theory game by a finite game and recall that optimal strategies of adversarial finite zero-sum games are mixed strategies \cite{Nash:1950, Nash:1951}.

Le Cam \cite{LeCam1955}, see also \cite{LeCam}, has substantially
extended  Wald's theory in the sense that  requirements of boundedness, or even finiteness,
of the loss function are replaced by  a requirement of lower semicontinuity, and the requirements  of
 the compactness of the decision space and the
  absolute continuity of the admissible measures with respect the Lebesgue measure are removed.
 These vast generalizations
come at some price of abstraction, yet reveal
 the relevance and utility of an appropriate complete Banach lattice of measures. In particular, this framework of Le Cam appears to facilitate efficient concrete
 approximation.

As an illustration,
let us describe a complete class theorem on a space of admissible measures, without the inclusion of
functions, where the observation map  consists of taking $n$-i.i.d.~samples, as in Equation \eqref{eq:stater0}.
Let $\mathcal{A} \subset \mathcal{M}(\X)$ be a subset of the Borel probability measures on a topological space
$\X$ and consider a quantity of interest $\Phi:\mathcal{A} \rightarrow \R$.
For $\mu \in \mathcal{A}$, the data $d$ is generated by i.i.d.~sampling with respect to
$\mu^{n}$. That is $d \sim \mu^{n}$. For $\mu^\dagger  \in \mathcal{A}$,
the statistical error $\Er\big(\theta, \mu^\dagger)$  of an estimator
  $\theta:\X^{n}\rightarrow \R$ of $\Phi(\mu^\dagger)$  is  defined in terms of a loss function
$V\colon\R\rightarrow \R$ as in \eqref{eq:stater0}. Define the least upper bound on that statistical error and the sharpest estimator as in \eqref{eq:stater1}
and \eqref{eq:stater2}.

Let $\Theta:=\{\theta:\X^{n}\rightarrow \R, \, \theta\,  \text{measurable}\}$ denote the space
of  estimators.
 Since,  in general
 the game
 $ \Er(\theta, \mu), \theta\in \Theta, \mu \in \mathcal{A}$  will not have a value, that is one will have
a  strict inequality
\[ \sup_{\mu \in \mathcal{A}}\inf_{\theta \in \Theta}{\Er(\theta, \mu) }<
 \inf_{\theta\in \Theta}\sup_{\mu \in \mathcal{A}}{\Er(\theta, \mu)}\, ,\]
 classical arguments in game theory suggest
that one extend to
 random estimators and random selection in $\mathcal{A}$. To that end,  let
the set of randomized estimators $\mathcal{R}:=\{\hat{\theta}:\X^{n}\times \mathcal{B}(\R) \rightarrow [0,1],\,  \hat{\theta}
\,\,  \text{Markov} \}$  be
 the set  of Markov kernels.
 To define a topology for $\mathcal{R}$,  define a  linear space of  measures as follows.
Let $\mathcal{A}^{n}:=\{\mu^{n}\in \mathcal{M}(\mathcal{X}^{n}):\mu \in \mathcal{A}\}$ denote the
corresponding set of measures generating sample data.
Say that $\mathcal{A}^{n}$ is dominated if  there exists  an $\omega\in \mathcal{M}(\mathcal{X}^{n})$ such that every $\mu^{n}\in \mathcal{A}^{n}$ is absolutely continuous with respect to $\omega$.
 According to the Halmos-Savage Lemma \cite{HalmosSavage},
 see also Strasser \cite[Lem.~20.3]{strasser1985mathematical}, the set $\mathcal{A}^{n}$ is dominated
if and only if there exists a countable mixture
$\mu^{*}:=\sum_{i=1}^{\infty}{\alpha_{i}\mu^{n}_{i}}$, with $\alpha_{i}\geq 0, \mu_{i}\in \mathcal{A}, i =1,\ldots \infty$,
and $\sum_{i=1}^{\infty}{\alpha_{i}}=1$, such that
$\mu^{n} \sim \mu^{*}, \mu \in \mathcal{A}$.
A  construct at the heart of Le Cam's approach is a natural linear space
notion of
a mixture space of $\mathcal{A}$, called
 the $L$-space
$L(\mathcal{A}^{n}):=L^{1}(\mu^{*})$. It follows easily,  see \cite[Lem.~41.1]{strasser1985mathematical},
 that $L(\mathcal{A}^{n})$ is the set of  signed measures which are absolutely continuous with respect
to $\mu^{*}$.  When $\mathcal{A}$ is not dominated a natural generalization of this construction
 \cite[Def.~41.3]{strasser1985mathematical} due to Le Cam \cite{LeCam1955} is used.  A crucial property
of the $L$-space  $L(\mathcal{A}^{n})$ is that not only is it a Banach lattice, see Strasser
 \cite[Cor.~41.4]{strasser1985mathematical}, but by \cite[Lem.~41.5]{strasser1985mathematical} it is a complete lattice.
The utility of the concept of a complete lattice to the complete class theorems can clearly be seen in
the proof of the lemma in Section 2 of   Wald and Wolfowitz' \cite{wald1951}
proof of  the complete class theorem when the number of decisions and the number of distributions is finite.
Then, the natural action of a randomized estimator on the bounded continuous function/mixture pairs
$C_{b}(R)\times L(\mathcal{A}^{n})$ is
\[f\hat{\theta}\nu :=\int{ \int {f(r)\hat{\theta}(x^{n},dr)}\nu(dx^{n})}, \quad f \in C_{b}(R),\, \nu \in L(\mathcal{A}^{n})\,   .\]
Let $\mathcal{R}$ be endowed with the topology of pointwise convergence with respect to this action, that is the
weak topology with respect to integration against $C_{b}(R)\times L(\mathcal{A}^{n})$.
Moreover, this weak topology also facilitates a definition of the space $\overline{\mathcal{R}}$
 of {\em generalized} random estimators
as bilinear  real-valued maps $\vartheta: C_{b}(R)\times L(\mathcal{A}^{n}) \rightarrow \R$ satisfying
$|\vartheta(f,\mu)| \leq \inorm{f}\snorm{\mu}$, $\vartheta(f,\mu)\geq 0$ for $f\geq 0$, $\mu \geq 0$,
and $\vartheta(1,\mu)=\mu(\X)$.
By \cite[Thm.~42.3]{strasser1985mathematical}
the set of generalized random estimators  $\overline{\mathcal{R}}$ is compact and convex
 and by \cite[Thm.~42.5]{strasser1985mathematical} of
Le Cam \cite{le1964sufficiency}, $\mathcal{R}$ is dense in $\overline{\mathcal{R}}$ in the weak
topology. Moreover, when $\mathcal{A}^{n}$ is dominated and one can restrict to a compact subset
$C\in \R$ of the decision space, then Strasser \cite[Cor.~42.8]{strasser1985mathematical} asserts
that $\overline{\mathcal{R}}=\mathcal{R}$.

Returning to our illustration,  if one let $W_{\mu},\mu \in \mathcal{A}$ be defined by
$W_{\mu}(r):=V(r-\Phi(\mu)), r \in \R, \mu \in \mathcal{A}$
 denote the associated family of loss functions, one can now define a generalization
of the statistical error function $\Er(\theta, \mu)$ of \eqref{eq:stater0}  to randomized
estimators $\hat{\theta}$ by
\[\Er(\hat{\theta}, \mu):= \int{ \int {W_{\mu}(r)\hat{\theta}(x^{n},dr)}\mu^{n}(dx^{n})},\quad  \hat{\theta} \in \mathcal{R}, \mu \in \mathcal{A}\, .\]
This definition reduces to the previous one \eqref{eq:stater0}  when the random estimator
$\hat{\theta}$ corresponds to a point estimator $\theta$ and extends naturally to $\overline{\mathcal{R}}$.
Finally, one says that an estimator $\vartheta^{*}\in \overline{\mathcal{R}}$ is Bayesian
if there exists a  probability  measure $m$ with finite support on $\mathcal{A}$ such that
\[ \int {\Er(\vartheta^{*}, \mu) m(d\mu)}   \leq \int {\Er(\vartheta, \mu) m(d\mu)}, \quad  \vartheta \in
\overline{\mathcal{R}}\, .  \]
The following complete class theorem follows from Strasser \cite[Thm.~47.9, Cor.~42.8]{strasser1985mathematical}
since one can naturally compactify the decision space $\R$ when the quantity of interest
$\Phi$ is bounded and the loss function $V$ is sublevel compact, that is has compact sublevel sets.
\begin{thm}
\label{thm_cc}
Suppose that the loss function $V$ is sublevel compact and
 the quantity of interest $\Phi:\mathcal{A} \rightarrow \R$ is bounded. Then,
 for each generalized randomized estimator $\vartheta \in \overline{\mathcal{R}}$, there exists a weak limit
$\vartheta^{*} \in \overline{\mathcal{R}}$ of Bayesian estimators
 such that
\[  \Er(\vartheta^{*}, \mu) \leq \Er(\vartheta, \mu), \quad \mu \in \mathcal{A}\, .\]
If, in addition, $\mathcal{A}$ is dominated, then there exists such a $\vartheta^{*} \in \mathcal{R}$.
\end{thm}
A comprehensive connection of these results, where Bayesian estimators
are defined only in terms of measures of finite support on $\mathcal{A}$,  with the framework of  Section \ref{sec_staterror_random} where
Bayesian estimators are defined in terms of Borel measures on
$\mathcal{A}$ is not available yet.
 Nevertheless it appears that much can be done in this regard. In particular, one can suspect that when
$\mathcal{A}$ is a  closed convex set  of probability measures equipped with the weak topology
 and $\X$ is a Borel subset of a Polish space, that  if
 the loss function $V$
is convex and $\Phi$ is affine and measurable,
 that the Choquet theory of Winkler  \cite{winkler1978integral,Winkler:1988}  can be used to
facilitate this connection.
Indeed, as mentioned above, complete class theorems are available for much more general loss functions
than continuous or convex,  more general decision spaces than $\R$,  and without absolute continuity assumptions.
 Moreover,  it is interesting
to note that, although randomization was introduced to obtain minmax results, when the loss function
$V$ is strictly convex, Bayesian estimators can be shown to be {\em  non-random}. This can be explicitly
observed in the definition  \eqref{eq:dsicees3repu} of Bayesian estimators when $V(x):=x^{2}$, and is understood
much more generally in Dvoretsky, Wald and Wolfowitz  \cite{DvoretzkyWaldWolfowitz}.
We conjecture that further simplifications can be obtained by allowing
{\em approximate} versions of complete class theorems, Bayesian estimators,  optimality,
and saddle points, as in Scovel, Hush and Steinwart's \cite{scovel2007approximate} extension
of classical Lagrangian duality theory to include approximations.

\section{Incorporating complexity and computation}

Although Decision Theory provides well posed notions of optimality and performance in statistical estimation,
it does not address the complexity of the actual computation of optimal or nearly optimal estimators and their evaluation against the data.
Indeed, although the abstract identification of an optimal estimator as the solution of an optimization problem provides a clear objective, practical applications require the actual implementation of the estimator on a machine
and its numerical evaluation against the data.

\subsection{Machine Wald}
The simultaneous  emphasis on {\em performance} and {\em computation} can be traced back to PAC (probably approximately correct) Learning initiated by Valiant \cite{Valiant} which has laid down the foundations of Machine Learning (ML). Indeed, as asserted by   Wasserman in his 2013 Lecture  ``The Rise of the Machines'' \cite[Sec.~1.5]{Wasserman_rise}:
\begin{quotation}
        \noindent
``There is another interesting difference that is worth pondering. Consider the problem of estimating a mixture of Gaussians. In Statistics we think of this as a solved problem. You use, for example, maximum likelihood which is implemented by the EM algorithm. But the EM algorithm does not solve the problem. There is no guarantee that the EM algorithm will actually find the MLE; it’s a shot in the dark. The same comment applies to MCMC methods.
In ML, when you say you’ve solved the problem, you mean that there is a polynomial time algorithm with provable guarantees. ''
\end{quotation}
  That is, on even par
with the rigorous performance analysis, Machine Learning also requires
that solutions be efficiently implementable on a computer, and often such efficiency is established by proving bounds on the amount of computation required to produce such a solution with a given algorithm.
Although Wald's theory of Optimal Statistical Decisions has  resulted in many important statistical discoveries, looking through the three
Lehmann symposia of Rojo and P{\'e}rez-Abreu \cite{Rojo:optimality1} in 2004,  and  Rojo \cite{Rojo:optimality2,Rojo:optimality3} in 2006 and 2009, it is clear that  the incorporation of the analysis  of  the computational algorithm, both in terms of its computational efficiency and its
 statistical optimality has not begun.
Therefore a natural answer to fundamental challenges in UQ appears to be the full incorporation of {\em computation} into a natural generalization of
 Wald's Statistical Decision Function framework, producing a framework one might call {\em Machine Wald}.

\subsection{Reduction Calculus}\label{redcalculus}
The resolution of minimax problems \eqref{eq:stater2scose} require, at an abstract level, searching in
the space of all possible functions of the data. By restricting models to the Bayesian class, the complete class theorem allows to limit this search to prior distributions on $\mathcal{A}$, i.e. to measures over spaces of measures and functions.
To enable the computation of these models  it is therefore necessary to identify conditions under which Minimax problems over measures over spaces of measures and functions can be reduced to the manipulation of finite-dimensional objects and develop the associated reduction calculus.
For  min or max problems over measures over spaces of measures (and possibly functions) this calculus can take the form of a reduction to a nesting of optimization problems over measures (and possibly functions for the inner part)  \cite{OwScSu2015, OwhadiScovel:2013, owhadiBayesiansirev2013}, which, in turn, can be reduced to searches over extreme points \cite{OSSMO:2011, Sullivan:2013, Hanconvexouq:2013, owhadiscovelexball2015}.

\subsection{Stopping conditions}
Many of these optimization problems will not be tractable. However even in the tractable case, which has rigorous
guarantees on the amount of computation required to obtain an approximate optima, it will be useful
to have stopping criteria for the specific algorithm and the specific problem instance under consideration, which can be used to guarantee when an approximate optima has been achieved. Although in the intractable case
no such guarantee will exist in general, intelligent choices of algorithms may result
in the attainment of approximate optima and such tests guarantee that fact.
 Ermoliev, Gaivoronski and Nedeva \cite{Ermoliev} successfully develop
such stopping criteria using Lagrangian  duality  and
 generalized Bender's decompositions by Geoffrion \cite{Geoffrion:Benders} for certain
Stochastic Optimization problems which are also relevant here.
In addition, the approximation of intractable problems by tractable ones will be important. Recently,
 Hanasusanto, Roitch, Kuhn  and Wiesemann \cite{hanasusanto2015distributionally}
 derive explicit conic reformulations for tractable problem classes and
suggest efficiently computable conservative approximations for intractable ones.

\subsection{On the Borel-Kolmogorov paradox}\label{subsecborelprad}
An oftentimes overlooked difficulty with Bayesian estimators lies in the fact that
for a prior $\pi \in \mathcal{M}(\mathcal{A})$, the posterior \eqref{eqcondexp} is not a measurable function of $\Ddata$ but a convex set $\Theta(\pi)$ of measurable functions $\theta$ of $\Ddata$ that are almost-surely equal to each other under the measure $\pi \cdot \Dmap$ on $\Dspace$.

A notorious pathological consequence is the Borel-Kolmogorov paradox (see Chapter 5 of \cite{Kolmogorov:2003} and Section 15.7 of \cite{Jaynes:2003}).
Recall that in the formulation of this paradox one considers the uniform distribution on the
 two-dimensional sphere and one is interested in obtaining the conditional distribution associated with a great circle of that sphere. If the problem is parameterized in spherical coordinates, then the resulting conditional distribution is uniform for the equator but non uniform for the longitude corresponding to the prime meridian.
The following theorem  suggests that this paradox is generic and dissipates the idea that it could be limited to fabricated toy examples. See also Singpurwalla and Swift \cite{SingpurwallaSwift:2002} for implications of this paradox in modeling and inference.

Recall that for $\pi \in \mathcal{M}(\mathcal{A})$, that $\Theta(\pi)$ is defined
as the convex set of measurable functions which are equal $\pi \cdot \Dmap$-everywhere to
 the regular conditional expectation \eqref{eqcondexp}.  Despite this indeterminateness
it is comforting to know that
\[\Er(\theta_2, \pi)=\Er(\theta_1, \pi), \quad \theta_{1},\theta_{2} \in \Theta(\pi)\, .\]
Moreover, it is also easy to see that if
 $\pi^\dagger$ is absolutely continuous with respect to $\pi$,
then
 $\theta_1(d)=\theta_2(d)$ with $\pi^\dagger \cdot \Dmap$  probability one
for all $\theta_1,\theta_2 \in \Theta(\pi)$, and consequently
\begin{equation*}\label{ded}
\Er(\theta_2, \pi^\dagger)=\Er(\theta_1, \pi^\dagger),\quad \theta_{1},\theta_{2} \in \Theta(\pi),\quad \pi^{\dagger}
 \prec \pi\,,
\end{equation*}
where the notation $\pi^{\dagger}
 \prec \pi$ means that $\pi^{\dagger}$ is absolutely continuous with respect to
 $\pi$.
The following theorem shows that this requirement of absolute continuity is necessary for all versions of conditional expectations $\theta \in \Theta(\pi)$ to share the same risk. See Subsection \ref{subsecborel}  for its proof.

\begin{thm}\label{thm:bidoneprior}
Assume that $V(x)=x^2$ and that the image  $\Phi(\mathcal{A})$ is a nontrivial interval.
If $\pi^\dagger$ is not absolutely continuous with respect to $\pi$ then
\begin{equation}\label{eqgjhghhjggyh}
\frac{1}{4}\leq \frac{\sup_{\theta_1,\theta_2 \in \Theta(\pi)} \big(\Er(\theta_2, \pi^\dagger)-\Er(\theta_1, \pi^\dagger)\big)}{ \big(\mathcal{U}(\mathcal{A})-\mathcal{L}(\mathcal{A})\big)^2 \sup_{B \in \mathcal{B}(\Dspace)\,:\, (\pi \cdot \Dmap)[B]=0} (\pi^\dagger \cdot \Dmap)[B]}  \leq 1
\end{equation}
where $\mathcal{U}(\mathcal{A})$ and $\mathcal{L}(\mathcal{A})$ are defined by \eqref{eq:defma1} and \eqref{eq:defma2}.
\end{thm}
\begin{rmk}\label{rmk:alwaysboundseed}
If moreover  $\pi^\dagger\cdot \Dmap$ is orthogonal to $\pi\cdot \Dmap$, that is there exists a set $B\in \mathcal{B}(\Dspace)$ such that $(\pi\cdot \Dmap)[B]=0$ and $(\pi^\dagger \cdot \Dmap)[B]=1$, then Theorem \ref{thm:bidoneprior} implies that
$\sup_{\theta_1,\theta_2 \in \Theta(\pi)}\big(\Er(\theta_2, \pi^\dagger)-\Er(\theta_1, \pi^\dagger)\big)$ is larger than  the statistical error of the midpoint estimator
\begin{equation*}\label{eq:midpoint}
\theta:=\frac{\mathcal{L}(\mathcal{A})+\mathcal{U}(\mathcal{A})}{2}.
\end{equation*}
\end{rmk}
As a remedy one can try (see \cite{Tjur:1974}, \cite{Tjur:1980} and \cite{Pfanzagl:1979}) constructing conditional expectations as disintegration or derivation limits defined as
\begin{equation}\label{eq:dsicrep0weq}
\E_{\pi \odot \Dmap}\big[\Phi(f,\mu)\big| \Drv=\Ddata\big]=\lim_{B\downarrow \{\Ddata\}} \E_{\pi \odot \Dmap}\big[\Phi(f,\mu)\big| \Drv\in B\big]
\end{equation}
where the limit $B\downarrow \{\Ddata\}$ is taken over a net of open neighborhoods of $\Ddata$. But as shown in \cite{KacSlepian:1959}, the limit generally depends on the net $B\downarrow \{\Ddata\}$ and the resulting conditional expectations can be distinctly different for different nets.
Furthermore the limit \eqref{eq:dsicrep0weq} may exist/not exist on subsets of $\Dspace$ of $(\pi\cdot\Dmap)$-measure  zero (which, as shown above, can lead to the inconsistency of the estimator). A related important issue is that  conditional probabilities can in general not be computed  \cite{AckFreerRoy:2010}. Observe that if the limit \eqref{eq:dsicrep0weq} does not exist then Bayesian estimation of $\Phi(f,\mu)$ may have significant oscillations as the precise measurement of $\Ddata$ becomes sharper.

\subsection{On Bayesian Robustness/Brittleness}

As much as classical numerical analysis shows that there are stable and unstable ways to discretize a partial differential equation,  positive \cite{Bernstein:1964, CastilloNickl:2013, Doob:1949, Kleijnvanderv:2012, LeCam:1953, Stuart:2010, vonMises:1964} and negative results \cite{Belot:2013, DiaconisFreedman:1986, Freedman:1963, Freedman:1999, Johnstone:2010, Leahu:2011, OwScSu2015, OwhadiScovel:2013, owhadiBayesiansirev2013,OwhadiScovel:2014} are forming an emerging understanding of stable and unstable ways to apply Bayes' rule in practice.
One aspect of stability concerns the sensitivity of  posterior conclusions with respect to the underlying models and prior beliefs.

\begin{quotation}
	\noindent ``Most statisticians would acknowledge that an analysis is not complete unless the sensitivity of the conclusions to the assumptions is investigated.  Yet, in practice, such sensitivity analyses are rarely used.  This is because sensitivity analyses involve difficult computations that must often be tailored to the specific problem.  This is especially true in Bayesian inference where the computations are already quite difficult.'' \cite{WassermanEtAl:1993}
\end{quotation}
Another aspect concerns situations where Bayes' rule is applied iteratively and posterior values become prior values for the next iteration. Observe in particular that when  posterior distributions (which are later on used as prior distributions) are only approximated (e.g.\ via MCMC methods), stability requires the convergence of the MCMC method in the same metric used to quantify the sensitivity of posterior with respect to the prior distributions.

In the context of the framework being developed here,
recent results \cite{OwScSu2015, OwhadiScovel:2013, owhadiBayesiansirev2013,OwhadiScovel:2014} on the extreme sensitivity (brittleness) of Bayesian inference in the TV and Prokhorov metrics appear to suggest that robust inference, in a continuous world under finite-information, should perhaps be done with reduced/coarse models rather than highly sophisticated/complex models (with a level of coarseness/reduction depending on the available finite-information) \cite{owhadiBayesiansirev2013}.

\subsection{Information Based Complexity}
From the point of view of practical applications it is clear that the set of possible models entering in the Minimax Problem \ref{eq:stater2scose} must be restricted by introducing constraints on computational complexity. For example, finding optimal models of materials in extreme environments is not the correct objective when these models  require full Quantum Mechanics calculations. A more productive approach is to search for computationally tractable optimal  models  in a given complexity class. Here one may wonder if Bayesian models remain a complete class for the resulting complexity constrained Minimax problems.
 It is also clear that computationally tractable optimal  models may not use all the available information, for instance a material model of bounded complexity should not use the state of every atom.
 The idea that fast computation requires computation with partial information forms the core of Information Based Complexity, the branch of computational complexity that studies the complexity of approximating continuous mathematical operations with discrete and finite ones up a to specified level of accuracy \cite{Traub1988, Woniakowski1986, Packel1987, Nemirovsky1992, Woniakowski2009}, where it is also augmented by concepts of contaminated and priced information associated with, for example, truncation errors and the cost of numerical operations.
Recent results \cite{OwhadiMultigrid:2015} suggest that Decision Theory concepts could be used, not only  to identify reduced models but also
 algorithms of near optimal complexity by reformulating the process of computing with partial information and limited  resources as that of playing underlying hierarchies of adversarial information games.

\section{Conclusion}
Although Uncertainty Quantification is still in its formative stage, much like the state of probability theory before its rigorous formulation by Kolmogorov in the 1930s, it has the potential to have an impact on the process of scientific discovery that is similar to the advent of scientific computing.
Its emergence remains sustained by the persistent need to make critical decisions with partial information and limited resources.
There are many paths to its development,  but one such path appears to be the incorporation of notions of computation and complexity in a generalization of Wald's decision framework built on Von Neumann's theory of adversarial games.

\section{Appendix}
\subsection{Construction of $\pi \odot \Dmap$}\label{subsecpidmap}
The below construction works when
 $\mathcal{A} \subseteq \mathcal{G}\times \mathcal{M}(\mathcal{X})$
for some Polish subset $\mathcal{G} \subset \mathcal{F}(\X)$ and $\mathcal{X}$ is Polish.
Observe that since $\Dspace$ is metrizable, it follows from \cite[Thm.~15.13]{AliprantisBorder:2006}, that, for any $B \in \mathcal{B}(\Dspace)$, the evaluation $\nu \mapsto \nu(B)$, $\nu \in \mathcal{M}(\Dspace)$,  is measurable.  Consequently, the measurability of $\Dmap$ implies that the mapping
\[
	\widehat{\Dmap}\colon \mathcal{A} \times \mathcal{B}(\Dspace) \to R
\]
defined by
\[
	\widehat{\Dmap}\bigl( (f,\mu), B\bigr) := \Dmap(f,\mu)[B], \quad \text{for } (f,\mu) \in  \mathcal{A}, B \in
 \mathcal{B}(\Dspace)
\]
is a transition function in the sense that, for fixed $(f,\mu) \in \mathcal{A}$,  $\widehat{\Dmap}\bigl( (f,\mu), \quark \bigr)$  is a probability measure, and, for fixed $B \in \mathcal{B}(\Dspace)$, $\widehat{\Dmap}\bigl( \quark ,B\bigr)$ is Borel measurable.
Therefore, by \cite[Thm.~10.7.2]{Bogachev2}, any $\pi \in \mathcal{M}(\mathcal{A})$, defines a probability measure
\[
	\pi \odot \Dmap \in \mathcal{M}\bigl(\mathcal{B}(\mathcal{A}) \times \mathcal{B}(\mathcal{D})\bigr)
\]
through
\begin{equation}
	\label{eq:palm}
	\pi\odot \Dmap \big[ A \times B \big] := \E_{(f.\mu) \sim \pi} \big[ \one_{A}(f, \mu) \Dmap(f, \mu)[B] \big],\quad \text{for } A \in \mathcal{B}(\mathcal{A}), B \in \mathcal{B}(\Dspace) ,
\end{equation}
where $\one_{A}$ is the indicator function of the set $A$:
\[
	\one_{A}(f,\mu) :=
	\begin{cases}
		1, & \text{if $(f,\mu) \in A$,} \\
		0, & \text{if $(f,\mu) \notin A$.}
	\end{cases}
\]
It is easy to see that $\pi$ is the $\mathcal{A}$-marginal of $\pi\odot \Dmap$.  Moreover, when $\mathcal{X}$ is Polish, \cite[Thm.~15.15]{AliprantisBorder:2006} implies that $\mathcal{M}(\mathcal{X})$ is Polish, and when $\mathcal{G}$ is Polish it follows that $\mathcal{A} \subseteq \mathcal{G}\times \mathcal{M}(\mathcal{X})$
is second countable.  Consequently, since $\Dspace$ is Suslin and hence second countable, it follows from \cite[Prop.~4.1.7]{Dudley:2002} that
\[
	\mathcal{B}\bigl(\mathcal{A} \times \Dspace\bigr)=\mathcal{B}(\mathcal{A}) \times\mathcal{B}(\Dspace)
\]
and hence $\pi \odot \Dmap$ is a probability measure on $\mathcal{A} \times \Dspace$.  That is,
\[
	\pi\odot \Dmap \in \mathcal{M}(\mathcal{A} \times \Dspace).
\]

Henceforth  denote $\pi\cdot\Dmap$  the corresponding  Bayes' sampling distribution defined by the $\Dspace$-marginal of $\pi\odot \Dmap$,  and note that, by \eqref{eq:palm}, one has
\begin{equation*}
    \label{eq:cdotexp}
    \pi\cdot\Dmap[B]:=\E_{(f,\mu)\sim \pi}\big[\Dmap(f,\mu)[B]\big], \quad \text{for } B \in
\mathcal{B}(\Dspace) .
\end{equation*}

Since both $\Dspace$ and $\mathcal{A}$ are Suslin it follows that the product $\mathcal{A} \times \Dspace$ is Suslin. Consequently, \cite[Cor.~10.4.6]{Bogachev2} asserts that regular conditional probabilities exist for any sub-$\sigma$-algebra of $\mathcal{B}\bigl(\mathcal{A} \times \Dspace\bigr)$. In particular, the product theorem of \cite[Thm.~10.4.11]{Bogachev2} asserts that product regular conditional probabilities
\[
	\bigl(\pi\odot \Dmap\bigr)|_{d} \in \mathcal{M}(\mathcal{A}), \quad \text{for } d \in \Dspace
\]
exist and that they are $\pi\cdot \Dmap$-a.e.\ unique.

\subsection{Proof of Theorem \ref{thm:bidoneprior}}\label{subsecborel}
If $\pi^\dagger\cdot \Dmap$ is not absolutely continuous with respect to $\pi \cdot \Dmap$ then there exists $B \in \mathcal{B}(\Dspace)$ such that
$(\pi \cdot \Dmap)[B]=0$ and $(\pi^\dagger \cdot \Dmap)[B]>0$. Let $\theta \in \Theta(\pi)$. Define
\begin{equation}
\theta_{y}(\Ddata):=\theta(\Ddata)1_{B^c}(\Ddata)+y 1_{B}(\Ddata)
\end{equation}
Then it is easy to see that if $y$ is in the range of $\Phi$ then $\theta_y\in \Theta(\pi)$. Now observe that for $y,z\in \Image(\Phi)$,
\begin{equation*}\label{eq:dsicee3srsesewep0}
\Er(\theta_y,\pi^\dagger)-\Er(\theta_z,\pi^\dagger)=\E_{(f,\mu,\Ddata)\sim \pi^\dagger\odot \Dmap}\Bigg[1_B(\Ddata)\Big(V\big(y-\Phi(f,\mu)\big)-V\big(z-\Phi(f,\mu)\big)\Big) \Bigg]
\end{equation*}
Hence, for $V(x)=x^2$, it holds true that
\begin{equation*}\label{eq:dwsicee3rsewesp0}
\Er(\theta_y,\pi^\dagger)-\Er(\theta_z,\pi^\dagger)=\big[(y-\gamma)^2-(z-\gamma)^2\big] (\pi^\dagger \cdot \Dmap)[B]
\end{equation*}
with
\begin{equation*}\label{eq:dsiddeeedcee3rsewesp0}
\gamma:=\E_{\pi^\dagger\odot \Dmap}[\Phi| \Drv \in B]
\end{equation*}
which proves

\begin{equation*}\label{eq:diggvdij3o}
\begin{split}
\sup_{\theta_2 \in \Theta(\pi)} \Er(\theta_2,\pi^\dagger)&-\inf_{\theta_1 \in \Theta(\pi)} \Er(\theta_1,\pi^\dagger) \geq  \sup_{B \in \mathcal{B}(\Dspace)\,:\, (\pi \cdot \Dmap)[B]=0,\,y,z\in \Image(\Phi)}
\\& \Big[\big(y-\E_{\pi^\dagger\odot \Dmap}[\Phi| \Drv \in B]\big)^2-\big(z-\E_{\pi^\dagger\odot \Dmap}[\Phi| \Drv \in B]\big)^2\Big] (\pi^\dagger \cdot \Dmap)[B],
\end{split}
\end{equation*}
and,
\begin{equation*}\label{eq:lowboedebay}
\begin{split}
\sup_{\theta_2 \in \Theta(\pi)} \Er(\theta_2,\pi^\dagger)&-\inf_{\theta_1 \in \Theta(\pi)} \Er(\theta_1,\pi^\dagger) \leq  \big(\mathcal{U}(\mathcal{A})-\mathcal{L}(\mathcal{A})\big)^2\sup_{B \in \mathcal{B}(\Dspace)\,:\, (\pi \cdot \Dmap)[B]=0} (\pi^\dagger \cdot \Dmap)[B].
\end{split}
\end{equation*}
To obtain the right hand side of \eqref{eqgjhghhjggyh}  observe that (see for instance \cite[Sec.~5]{Doob:1994})  there exists  $B^*\in \mathcal{B}(\Dspace)$ such that
\begin{equation*}
(\pi^\dagger \cdot \Dmap)[B^*]=\sup_{B \in \mathcal{B}(\Dspace)\,:\, (\pi \cdot \Dmap)[B]=0} (\pi^\dagger \cdot \Dmap)[B]
\end{equation*}
and (since $\theta_2=\theta_1$ on the complement of $B^*$)
\begin{equation*}\label{eq:dsicssrsesewep0}
\begin{split}
\sup_{\theta_1, \theta_2 \in \Theta(\pi)} &\big(\Er(\theta_2,\pi^\dagger)-\Er(\theta_1,\pi^\dagger)\big)\\&=\sup_{\theta_1, \theta_2 \in \Theta(\pi)}\E_{(f,\mu,\Ddata)\sim \pi^\dagger\odot \Dmap}\Bigg[1_{B^*}(\Ddata)\Big(V\big(\theta_2-\Phi(f,\mu)\big)-V\big(\theta_1-\Phi(f,\mu)\big)\Big) \Bigg].
\end{split}
\end{equation*}
We conclude by observing that for $V(x)=x^2$,
\begin{equation*}
\sup_{\theta_1, \theta_2 \in \Theta(\pi)} \Big(V\big(\theta_2-\Phi(f,\mu)\big)-V\big(\theta_1-\Phi(f,\mu)\big)\Big)\leq \big(\mathcal{U}(\mathcal{A})-\mathcal{L}(\mathcal{A})\big)^2.
\end{equation*}

\subsection{Conditional expectation as an orthogonal projection}\label{subseccondexp}

It easily follows from Tonelli's Theorem
that
\begin{equation*}
\E_{\pi \cdot \Dmap}[h^{2}]=\E_{\pi \odot \Dmap}[h^{2}]=\E_{(f,\mu)\sim \pi}\E_{ \Dmap(f,\mu)}[h^{2}]\, .
\end{equation*}
By considering the sub $\s$-algebra $\mathcal{A} \times \mathcal{B}(\Dspace) \subset
\mathcal{B}(\mathcal{A}\times \Dspace) =\mathcal{B}(\mathcal{A}) \times \mathcal{B}(\Dspace)$
 it follows from \emph{e.g.}\ Theorem 10.2.9 of \cite{Dudley:2002} that
$L^2_{\pi\cdot \Dmap}(\Dspace)$ is a closed Hilbert subspace of the Hilbert space
 $L^2_{\pi \odot \Dmap}(\mathcal{A}\times \Dspace)$  and
 the conditional expectation of $\Phi$ given the random variable $ \Drv$
 is the orthogonal projection from $L^2_{\pi \odot \Dmap}(\mathcal{A}\times \Dspace)$ to $L^2_{\pi\cdot \Dmap}(\Dspace)$.

%\paragraph{Acknowledgments.} The authors gratefully acknowledge this work supported by  the Air Force Office of Scientific Research under Award Number  FA9550-12-1-0389 (Scientific Computation of Optimal Statistical Estimators) and the U.S. Department of Energy Office of Science, Office of Advanced Scientific Computing Research, through the Exascale Co-Design Center for Materials in Extreme Environments (ExMatEx, LANL Contract No DE-AC52-06NA25396, Caltech Subcontract Number 273448).

\bibliographystyle{plain}
\bibliography{merged}

\end{document}